\documentclass[12pt,a4paper]{article}
\usepackage{a4wide}
\usepackage{amsfonts,amsmath}
\usepackage{latexsym}
\usepackage{times}
\author{K\'aroly J. B\"or\"oczky, Oriol Serra}
\title{Remarks on the equality case of the Bonnesen inequality }

\newcommand{\proof}{\noindent{\it Proof: }}
\newcommand{\proofbox}{\mbox{ $\Box$}\\}
\newcommand{\R}{\mathbb{R}}

\newtheorem{lemma}{Lemma}[section]
\newtheorem{theo}[lemma]{Theorem}

\newtheorem{claim}[lemma]{Claim}

\begin{document}

\maketitle

\begin{abstract}
An argument is provided for the
 equality case of the high dimensional  Bonnesen inequality
for sections. The known equality case of the  Bonnesen inequality for projections
is presented as a consequence.

\end{abstract}

\section{Introduction}

We write $\mu_d$ for the $d$-dimensional Lebesgue measure.
Let $S^{d-1}$ be the unit sphere in $\R^d$.
For a linear subspace $\Pi$ of $\R^d$, the orthogonal projection into $\Pi$
is denoted by $p_\Pi$. In the special case when $\Pi=u^\bot$
for a $u\in S^{d-1}$, the orthogonal projection into $u^\bot$ is denoted by $\pi_u$.
In addition, the  convex hull of $x_1,\ldots,x_k$
is denoted by $[x_1,\ldots,x_k]$.

The results in this note belong to the very heart of the Brunn-Minkowski theory,
so any of the monographs T. Bonnesen, W. Fenchel \cite{BoF87},
P.M. Gruber \cite{Gru07} and R. Schneider \cite{Sch93},
or the survey paper R.J. Gardner \cite{Gar02} provide the sufficient background.

Let $A$ and $B$ be convex bodies (compact convex sets with non-empty interiors)
in $\R^d$ for this section.
The Brunn-Minkowski inequality states
\begin{theo}[Brunn-Minkowski]
If $\alpha,\beta>0$, then
\label{brunnmink}
$$
\mu_d(\alpha\,A+\beta\,B)\geq
\left(\alpha\,\mu_d(A)^{\frac1d}+\beta\,\mu_d(B)^{\frac1d}\right)^d,
$$
with equality if and only if $A$ and $B$ are homothetic.
\end{theo}

According to the H\"older inequality, if $M,N>0$, then
$$
 \left(\alpha\,M^{\frac1{d-1}}+\beta\,N^{\frac1{d-1}}\right)^{d-1}
 \left(\alpha\,\frac{\mu_d(A)}M+\beta\,\frac{\mu_d(B)}N\right)\geq
\left(\alpha\,\mu_d(A)^{\frac1d}+\beta\,\mu_d(B)^{\frac1d}\right)^d,
$$
with equality if and only if $\frac{\mu_d(A)^{\frac1d}}{M^{\frac1{d-1}}}=
\frac{\mu_d(B)^{\frac1d}}{N^{\frac1{d-1}}}$.
Therefore the following result due to T. Bonnesen \cite{Bon29}
 strengthens the Brunn-Minkowki inequality.

\begin{theo}[Bonnesen I]
\label{Bonnesen-section}
If for a linear $(d-1)$-space $L$ in $\R^d$,  $M$ and $N$ are the maximal $(d-1)$-volumes
of the sections of $A$ and $B$, respectively, by hyperplanes parallel to $L$, then
$$
\mu_d(\alpha\,A+\beta\,B)\geq
 \left(\alpha\,M^{\frac1{d-1}}+\beta\,N^{\frac1{d-1}}\right)^{d-1}
 \left(\alpha\,\frac{\mu_d(A)}M+\beta\,\frac{\mu_d(B)}N\right).
$$
\end{theo}

Theorem~\ref{Bonnesen-section} has the following consequence about projections
(see also Section~\ref{secproj}).

\begin{theo}[Bonnesen II]
\label{Bonnesen-proj}
For $u\in S^{d-1}$, if  $M=\mu_{d-1}(\pi_uA)$ and $N=\mu_{d-1}(\pi_uB)$, then
$$
\mu_d(\alpha\,A+\beta\,B)\geq
 \left(\alpha\,M^{\frac1{d-1}}+\beta\,N^{\frac1{d-1}}\right)^{d-1}
 \left(\alpha\,\frac{\mu_d(A)}M+\beta\,\frac{\mu_d(B)}N\right).
$$
\end{theo}

The goal of this note is to characterize the equality cases in
Bonnesen's inequalities Theorems~\ref{Bonnesen-section}
and \ref{Bonnesen-proj}.
We use the notations of these theorems.
We note that Theorem~\ref{Bonnesen-proj-equa},
and the two dimensional case of Theorem~\ref{Bonnesen-section-equa} are proved
by G.A. Freiman, D. Grynkiewicz, O. Serra, Y.V. Stanchescu \cite{FGSS}.

For $u\in S^{d-1}$,
we say that a convex body $K$ is obtained from a convex body $C$ by stretching along
$u$, if there exist $\lambda\geq 0$ and $w\in\R^d$ such that
$K=C+[w,w+\lambda u]$. In particular $K=C+w$ if $\lambda=0$.

\begin{theo}
\label{Bonnesen-section-equa}
Equality holds in Theorem~\ref{Bonnesen-section} if and only if either $A$ and $B$
are homothetic, or there exist $v\in S^{d-1}$, homothetic convex bodies $A'$ and $B'$,
and hyperplane $H$ parallel to $L$, such that $\pi_v(A')=\pi_v(A'\cap H)$,
and $A$ and $B$ are obtained from $A'$ and $B'$, respectively, by stretching
along $v$.
\end{theo}

We note that the condition $\pi_v(A')=\pi_v(A'\cap H)$ is equivalent saying
that $A'\subset (A'\cap H)+\R v$. Convex bodies for which there exist such
hyperplane $H$ and unit vector $v$ are characterized in M. Meyer \cite{Mey98}.

As we discuss in Section~\ref{secproj}, the following is a simple consequence
of Theorem~\ref{Bonnesen-section-equa} {\it via} Steiner symmetrization.

\begin{theo}[Freiman,Grynkiewicz,Serra,Stanchescu]
\label{Bonnesen-proj-equa}
Equality holds in Theorem~\ref{Bonnesen-proj} if and only if
there exist homothetic convex bodies $A'$ and $B'$ such that
 $A$ and $B$ are obtained from $A'$ and $B'$, respectively, by stretching along $u$.
\end{theo}

Our proofs of the two inequalities by Bonnesen, and the characterizations of
the equality cases are based on the $(d-1)$-dimensional Brunn-Minkowski inequality, and its equality case. Therefore we provide a new proof for the $d$-dimensional
Brunn-Minkowski inequality and its equality case.

As related results, a true  discrete analogue of the Bonnesen inequality in the plane is proved by
D. Grynkiewicz, O. Serra \cite{GrS10}, and the equality conditions are clarified by
G. A. Freiman, D. Grynkiewicz, O. Serra, Y. V. Stanchescu \cite{FGSSI}.
In addition, M. Meyer \cite{Mey98} proves a crucial property of  a given convex body's sections of maximal $(d-1)$-volume parallel
to a hyperplane.

\section{Minkowski linear combinations}
\label{secsum}

In this section we recall some well-known simple but useful observations
about Minkowski linear combinations of convex bodies
(see P.M. Gruber \cite{Gru07} or R. Schneider \cite{Sch93}). If $X$ is a compact convex set
in $\R^d$, then its support function is
$$
h_X(v)=\max_{x\in X}\langle v,x\rangle \mbox{ \ for $v\in\R^d$}.
$$
Then $h_X$ is a positive homogeneous and convex function on $\R^d$,
which determines $X$ uniquely. In addition, if $Y$ is another compact convex set,
$\Pi$ is a linear subspace, and $\alpha,\beta>0$, then
\begin{eqnarray}
\label{supportsum}
h_{\alpha X+\beta Y}&=&\alpha\, h_X+\beta\, h_Y\\
\label{projsum}
p_\Pi(\alpha X+\beta Y)&=&\alpha\, p_\Pi X+\beta\, p_\Pi Y.
\end{eqnarray}
We note that
$v\in S^{d-1}$ is  exterior unit normal vector to a convex body $K$ in $\R^d$
at $x\in K$ if and only if $\langle v,x\rangle=h_K(v)$.  The following
is a simple but useful consequence of (\ref{supportsum}).

\begin{claim}
\label{Minkowskisum}
Let $C=\alpha\,A+\beta\,B$ for convex bodies $A,B$ in $\R^d$ and $\alpha,\beta>0$,
and let $z_0=\alpha\,x_0+\beta\,y_0$ for $z_0\in C$, $x_0\in A$, $y_0\in B$.
\begin{description}
\item{(i)}  If $x_0\in\partial A$ and $y_0\in\partial B$ with exterior unit normal vector $v$,
then $z_0\in\partial C$ with exterior unit normal vector $v$.
\item{(ii)} If $z_0\in\partial C$ with exterior unit normal vector $v$, then $x_0\in\partial A$ and
$y_0\in\partial B$ with exterior unit normal vector $v$.
\end{description}
\end{claim}

Our first application of Claim~\ref{Minkowskisum} is about planar convex bodies.

\begin{claim}
\label{planarsum}
Let $l$ be a line in $\R^2$ with $0\in l$, and let $C=\alpha\,A+\beta\,B$ for
convex bodies $A,B$ in $\R^2$ and $\alpha,\beta>0$
with $\alpha+\beta=1$. In addition, we assume that
\begin{description}
\item{(i)}  for any $z\in C$, $z+l$ intersects $A$ and $B$, and
$C\cap [z+l]=\alpha(A\cap [z+l])+\beta(B\cap [z+l])$,
\item{(ii)} there exists $z\in C$ such that
$C\cap [z+l]=A\cap [z+l]=B\cap [z+l]$.
\end{description}
Then $A=B$.
\end{claim}
\proof Let $l=\R u$ for the unit vector $u$, and let $v\in u^\bot$ be a unit vector.
In this case $\pi_uA=\pi_uB=\pi_uC=[av,bv]$ for some $a<b$.
There exist convex functions $f,g,\varphi,\psi$ on $[a,b]$ such that
\begin{eqnarray*}
A&=&\{tv+su:\,a\leq t\leq b\mbox{ and }-g(t)\leq s\leq f(t)\}\\
B&=&\{tv+su:\,a\leq t\leq b\mbox{ and }-\psi(t)\leq s\leq \varphi(t)\}.
\end{eqnarray*}
It follows from condition (i) and from Claim~\ref{Minkowskisum} that
$f'(t)=\varphi'(t)$ and $g'(t)=\psi'(t)$ wherever the derivatives exist,
thus there exist constants $\gamma,\delta$ such that
$f(t)=\varphi(t)+\gamma$ and $g(t)=\psi(t)+\delta$ for $t\in[a,b]$.
However condition (ii) yields that $\gamma=\delta=0$, therefore $A=B$.
\proofbox

Let $K$ be a convex body in $\R^d$, and let $u\in S^{d-1}$. For each line $l$ parallel with $u$
and intersecting ${\rm int}K$, we translate the segment $l\cap K$ along $l$ into the position
where the midpoint of the translated segment lies in $u^\bot$. The closure of the  union
of these translated segments is the Steiner symmetrial ${\cal S}_uK$ of $K$.
For another representation of the Steiner symmetrization,
we note that there exist concave functions $f$ and $g$
on $\pi_u(K)$ such that
$$
K=\{x+\lambda u:\,x\in\pi_uK\mbox{ and }-g(x)\leq \lambda\leq f(x)\}.
$$
Then
\begin{equation}
\label{steiner-func}
{\cal S}_uK=\left\{x+\lambda u:\,x\in\pi_u(K)\mbox{ and }
\frac{-f(x)-g(x)}2\leq \lambda\leq \frac{f(x)+g(x)}2\right\}.
\end{equation}
It follows that ${\cal S}_uK$ is a convex body symmetric through $u^\bot$, and
$\mu_d({\cal S}_uK)=\mu_d(K)$.

\begin{claim}
\label{steiner-symm}
For convex bodies $A$ and $B$ in $\R^d$, $u\in S^{d-1}$,
and $\alpha,\beta>0$, we have
$$
\alpha\,{\cal S}_uA+\beta\,{\cal S}_uB\subset{\cal S}_u(\alpha\,A+\beta\,B).
$$
In addition, if equality holds, and $a$ and $b$ are lines parallel to $u$ intersecting
${\rm int}\,A$ and ${\rm int}\,B$, and  there exist
parallel supporting hyperplanes at the top endpoints of
$a\cap {\cal S}_uA$ and $b\cap {\cal S}_uB$ to ${\cal S}_uA$ and
to ${\cal S}_uB$, respectively, then
there exist
parallel supporting hyperplanes at the top endpoints of
$a\cap A$ and $b\cap B$,
and parallel supporting hyperplanes at the bottom endpoints of
$a\cap A$ and $b\cap B$ to $A$ and
to $B$, respectively.
\end{claim}
\proof Let $l$ be a line parallel with $u$
and intersecting ${\rm int}(\alpha\,A+\beta\,B)$, and let $z_0$ be one of the endpoints of
$l\cap(\alpha\,{\cal S}_uA+\beta\,{\cal S}_uB)$. It follows by
Claim~\ref{Minkowskisum} (ii) that  $z_0=\alpha\,x_0+\beta\,y_0$, where
$x_0$ and $y_0$ are boundary points of ${\cal S}_uA$ and ${\cal S}_uB$,
sharing a common exterior unit vector with $z_0$. Therefore
$a=x_0+\R u$ and $b=y_0+\R u$ satisfy $l=\alpha\,a+\beta\,b$ and
$$
l\cap(\alpha\,{\cal S}_uA+\beta\,{\cal S}_uB)=
\alpha\,(a\cap{\cal S}_uA)+\beta\,(b\cap {\cal S}_uB).
$$
In particular
\begin{eqnarray}
\nonumber
\mu_1\left(l\cap{\cal S}_u(\alpha\,A+\beta\,B)\right)&=&
\mu_1\left(l\cap(\alpha\,A+\beta\,B)\right)\\
\label{lineq}
&\geq&\alpha\,\mu_1(a\cap A)+\beta\,\mu_1(b\cap B)\\
\nonumber
&=&
\alpha\,\mu_1(a\cap{\cal S}_uA)+\beta\,\mu_1(b\cap {\cal S}_uB)\\
\nonumber
&=&\mu_1\left(l\cap[\alpha\,{\cal S}_uA+\beta\,{\cal S}_uB]\right),
\end{eqnarray}
which in turn yields
$\alpha\,{\cal S}_uA+\beta\,{\cal S}_uB\subset{\cal S}_u(\alpha\,A+\beta\,B)$.

Assume now that $\alpha\,{\cal S}_uA+\beta\,{\cal S}_uB={\cal S}_u(\alpha\,A+\beta\,B)$,
and hence equality holds in (\ref{lineq}) for any  line $l$ parallel with $u$
and intersecting ${\rm int}(\alpha\,A+\beta\,B)$. It follows  that
\begin{equation}
\label{lAB}
l\cap(\alpha\,A+\beta\,B)=\alpha(a\cap A)+\beta(b\cap B).
\end{equation}
Writing $x_1,y_1,z_1$ to denote the top endpoint, and 
$x_2,y_2,z_2$ to denote the bottom endpoint of
 $a\cap A$, $b\cap B$ and $l\cap(\alpha\,A+\beta\,B)$, we deduce $z_i=\alpha\,x_i+\beta\,y_i$ for $i=1,2$, from (\ref{lAB}). Therefore
 Claim~\ref{Minkowskisum} (ii) completes the argument.
\proofbox

To introduce another method of symmetrization,
 let $K$ be a convex body in $\R^d$, and let $l$ be a line. For each
hyperplane $H$ orthogonal to $l$ and intersecting ${\rm int}K$,
consider the $(d-1)$-ball in $H$ with the same $(d-1)$-volume
as $H\cap K$ and centred at $H\cap l$. The closure of the union of these
$(d-1)$-balls centred on $l$ is a convex body ${\cal R}_lK$
by the Brunn-Minkowski inequality,
and ${\cal R}_lK$ is called the Schwarz-rounding of $K$. Readily
$\mu_d({\cal R}_lK)=\mu_d(K)$.
Similar argument to the one for Claim~\ref{steiner-symm}
(or using the fact that the Schwarz-rounding can be obtained
as the limit of repeated Steiner symmetrizations through
hyperplanes containing $l$) yields

\begin{claim}
\label{schwarz-round}
For convex bodies $A$ and $B$ in $\R^d$, line $l$,
and $\alpha,\beta>0$, we have
$$
\alpha\,{\cal R}_lA+\beta\,{\cal R}_lB\subset{\cal R}_l(\alpha\,A+\beta\,B).
$$
\end{claim}

Schwarz-rounding will be a basic tool for our proof of Theorem~\ref{Bonnesen-section-equa}.
It was W. Blaschke who gave a simple proof of the Brunn-Minkowski inequality using Schwarz rounding in \cite{Bla16}.

\section{Proof of Theorem~\ref{Bonnesen-section-equa}}

If the conditions stated in Theorem~\ref{Bonnesen-section-equa} hold, then
we readily have equality in Theorem~\ref{Bonnesen-section}.
For the reverse statement, we subdivide the argument into three sections.

\subsection{ A little preparation}

First we introduce some notation.
Let $u\in S^{d-1}$ be orthogonal to $L$,
let $K$ be a convex body in $\R^d$ , and let $Q$ be the
 the maximal $(d-1)$-volume
of the sections of $K$ by hyperplanes parallel to $L$.
For $s\in(0,Q)$, let
\begin{eqnarray*}
k_-(s)&=&\min\left\{p:\,K\cap (pu+L)\neq\emptyset
\mbox{ and } \mu_{d-1}(K\cap (pu+L))\geq s\right\} \\
k_+(s)&=&\max\left\{p:\,K\cap (pu+L)\neq\emptyset
\mbox{ and } \mu_{d-1}(K\cap (pu+L))\geq s\right\}.
\end{eqnarray*}
In addition we define
$$
\mbox{$K_-(s)=K\cap (k_-(s)u+L)$ and
$K_+(s)=K\cap (k_+(s)u+L)$}.
$$
We observe that $K\cap (pu+L)\neq\emptyset$
 if and only if $p\in[k_-(0),k_+(0)]$, possibly $k_-(Q)=k_+(Q)$,
but $k_-(s)<k_+(s)$ if $s<Q$.
It follows from the $(d-1)$-dimensional case of the Brunn-Minkowski inequality that
\begin{equation}
\label{Ksects}
\mu_{d-1}(K\cap (pu+L))\geq s
\mbox{ \ for $s\in(0,Q]$ if and only if $p\in[k_-(s),k_+(s)]$}.
\end{equation}
We observe that if the "top"and "bottom" sections of $K$ parallel to $L$ are of zero 
 $\mu_{d-1}$-measure, then $\mu_{d-1}(K_+(s))=\mu_{d-1}(K_-(s))=s$ for 
$s\in[0,Q]$. In general, we have
\begin{eqnarray}
\label{K-slarge}
\mbox{if \ }\mu_{d-1}(K_-(s))&>&s,
\mbox{ \ then $k_-(s)=k_-(0)$}  \\
\label{K+slarge}
\mbox{if \ }\mu_{d-1}(K_+(s))&>&s,
\mbox{ \ then $k_+(s)=k_+(0)$}.
\end{eqnarray}

Calculating the integral of $f(p)=\mu_{d-1}(K\cap (pu+L))$ for
$p\in[k_-(0),k_+(0)]$ by calculating the area of the part of $\R^2$ between the
graph of $f$  and the first axis using Fubini's theorem, and after that  using (\ref{Ksects}) yield
\begin{eqnarray}
\nonumber
\mu_d(K)&=&\int_{k_-(0)}^{k_+(0)}\mu_{d-1}(K\cap (pu+L))\,dp=
\int_0^Q\mu_1(\{p\in \R:\,\mu_{d-1}(K\cap (pu+L))\geq s\})\,ds\\
\label{V(K)}
&=&\int_0^Q\left(k_+(s)-k_-(s)\right)\,ds.
\end{eqnarray}

As the final part of our preparation, we discuss the case
$k_-(Q)<k_+(Q)$. We have equality in (\ref{Ksects})
for $s=Q$, therefore the equality case of the Brunn-Minkowski inequality implies
that $K\cap (pu+L)$ is a translate of $K_-(Q)$ for
$p\in[k_-(Q),k_+(Q)]$. Let
$K_+(Q)=K_-(Q)+\lambda v$ for $v\in S^{d-1}$
and $\lambda>0$. It follows  by the convexity of $K$
that $k_-(Q)<k_+(Q)$ implies
\begin{equation}
\label{Kstretch0}
\{x\in K:\,[k_-(Q)\leq\langle x,u\rangle\leq k_+(Q)\}=K_-(Q)+[0,\lambda v]
\mbox{ \ and $K\subset K_-(Q)+\R v$},
\end{equation}
and in turn
\begin{equation}
\label{Kstretch1}
\pi_v(K)=\pi_v(K_-(Q)),
\end{equation}
and that $K$ is obtained from the convex body
\begin{equation}
\label{Kstretch2}
K'=\bigcup_{s\in[0,Q]}\left((K_+(s)-\lambda v)\cup K_-(s)\right)
\end{equation}
by stretching along $v$.

\subsection{ A proof of Theorem~\ref{Bonnesen-section}}

Replacing $A$ and $B$ by $M^{\frac{-1}{d-1}}\,A$ and $N^{\frac{-1}{d-1}}\,B$,
if necessary, we may assume that
\begin{equation}
\label{MN1}
M=N=1.
\end{equation}
Let $C=\alpha A+\beta B$, and we write
$a_-(s),a_+(s),A_-(s),A_+(s)$, or $b_-(s),b_+(s),B_-(s),B_+(s)$,
or $c_-(s),c_+(s),C_-(s),C_+(s)$ to denote $k_-(s),k_+(s),K_-(s),K_+(s)$
if $K=A$, or $K=B$, or $K=C$, respectively.
We observe that
if $t\in(0,1]$, then
$$
\alpha A_+(t)+\beta B_+(t)\subset C\cap([\alpha a_+(t)+\beta b_+(t)]\,u+L].
$$
Therefore (\ref{Ksects}), the analogous relation for
$A_-(t)$ and $B_-(t)$,
 and the $(d-1)$-dimensional case of the Brunn-Minkowski inequality
 yield that
\begin{eqnarray}
\label{cab+}
c_+([\alpha+\beta]^{d-1}t)&\geq& \alpha a_+(t)+\beta b_+(t)\\
\label{cab-}
c_-([\alpha+\beta]^{d-1}t)&\leq& \alpha a_-(t)+\beta b_-(t).
\end{eqnarray}
We deduce by (\ref{V(K)}) that
\begin{eqnarray}
\label{muC1}
\mu_d(C)&\geq&\int_0^{(\alpha+\beta)^{d-1}}
\left(c_+(s)-c_-(s)\right)\,ds\\
\label{muC2}
&\geq&(\alpha+\beta)^{d-1}\int_0^1[\alpha a_+(t)+\beta b_+(t)]-
[ \alpha a_-(t)+\beta b_-(t)]\,dt\\
\label{muC}
&=&(\alpha+\beta)^{d-1}\left[\alpha\,\mu_d(A)+\beta\,\mu_d(B)\right] .
\end{eqnarray}

\subsection{Analyzing the equality case}

To simplify the formulae, in addition to  (\ref{MN1}), we also assume
\begin{equation}
\label{alphabeta1}
\alpha+\beta=1.
\end{equation}
Let us assume that
\begin{equation}
\label{bonnesen-sect-equa}
\mu_d(C)=\alpha\,\mu_d(A)+\beta\,\mu_d(B),
\end{equation}
and hence equality holds in (\ref{cab+}) and (\ref{cab-})
for $t\in(0,1]$.
In particular, there exists no $p>\alpha a_+(t)+\beta b_+(t)$
such that
$$
\mu_{d-1}(C\cap[pu+L])\geq t.
$$
Using the $(d-1)$-dimensional the Brunn-Minkowski inequality
and its equality case,
and that $\mu_{d-1}(A_+(t))=t=\mu_{d-1}(B_+(t))$
if $a_+(t)<a_+(0)$ and $b_+(t)<b_+(0)$, we deduce
\begin{eqnarray}
\label{CAB+}
C_+(t)&=&\alpha\,A_+(t)+\beta\, B_+(t),
\mbox{ \ and $A_+(t)$ and $B_+(t)$ are translates},\\
\label{CAB-}
C_-(t)&=&\alpha\,A_-(t)+\beta\, B_-(t),
\mbox{ \ and $A_-(t)$ and $B_-(t)$ are translates}
\end{eqnarray}
for all $t\in(0,1]$. Thus we may assume that
\begin{equation}
\label{AB1}
A_-(1)=B_-(1)\subset u^\bot.
\end{equation}

We note that equality holds in (\ref{muC1}), as well, therefore
\begin{equation}
\label{Cmax}
\mbox{$C_+(1)$ and $C_-(1)$ are sections of $C$ of maximal $(d-1)$-volume
among the ones parallel to $L$.}
\end{equation}
For the final part of the argument, we distinguish cases depending on
whether the section of maximal $(d-1)$-volume is unique.\\

\noindent{\bf Case 1 } $a_+(1)=a_-(1)$ and $b_+(1)=b_-(1)$.

First we show that
\begin{equation}
\label{abt}
\mbox{$a_+(t)=b_+(t)$ and $a_-(t)=b_-(t)$ for $t\in[0,1]$}.
\end{equation}
We observe that $a_+(0)=a_+(1)$ is equivalent saying
that the top section of $A$ parallel to $L$
is a section of maximal $(d-1)$-volume.
Possibly after reversing $u$, we may assume that $a_+(0)>a_+(1)$.
Let $t_+\in[0,1)$ be the maximal $t\in[0,1)$
such that $a_+(t)=a_+(0)$,
and let $t_-\in[0,1]$ be the maximal $t\in[0,1]$
such that $a_-(t)=a_-(0)$.

Let $\widetilde{A}$, $\widetilde{B}$ and $\widetilde{C}$
be the Schwarz rounding of $A$, $B$ and $C$ with respect to
$\R u$. In particular, (\ref{MN1}) yields that the maximal $(d-1)$-volumes
of the sections of
$\widetilde{A}$ and $\widetilde{B}$ parallel to $L$ are $1$.
 It follows from the Bonnesen inequality (\ref{muC}),
from the assumption of equality (\ref{bonnesen-sect-equa}), and
Claim~\ref{schwarz-round} that
\begin{eqnarray*}
\alpha\,\mu_d(A)+\beta\,\mu_d(B)&=&
 \mu_d(C)=
\mu_d(\widetilde{C})\geq \mu_d(\alpha\,\widetilde{A}+\beta\,\widetilde{B})\\
&\geq&
\alpha\,\mu_d(\widetilde{A})+\beta\,\mu_d(\widetilde{B})=
\alpha\,\mu_d(A)+\beta\,\mu_d(B).
\end{eqnarray*}
Therefore $\mu_d(\widetilde{C})= \mu_d(\alpha\,\widetilde{A}+\beta\,\widetilde{B})$,
and hence Claim~\ref{schwarz-round} yields
$$
\widetilde{C}=\alpha\,\widetilde{A}+\beta\,\widetilde{B}.
$$
We define
$\tilde{a}_-(t)$, $\tilde{a}_+(t)$, $\widetilde{A}_-(t)$,
$\widetilde{A}_+(t)$, or
$\tilde{b}_-(t)$, $\tilde{b}_+(t)$, $\widetilde{B}_-(t)$, $\widetilde{B}_+(t)$,
or $\tilde{c}_-(t)$, $\tilde{c}_+(t)$, $\widetilde{C}_-(t)$, $\widetilde{C}_+(t)$ to denote $k_-(t),k_+(t),K_-(t),K_+(t)$
if $K=\widetilde{A}$, or $K=\widetilde{B}$, or $K=\widetilde{C}$, respectively.
We observe that $\tilde{a}_+(t)=a_+(t)$ for $t\in[0,1]$, and $\widetilde{A}_+(t)$
is a $(d-1)$-ball with $\mu_{d-1}(\widetilde{A}_+(t))=t$
for $t\in[t_+,1]$, and we have the
similar statements for the analogous quantities.
Since  $\mu_d(\widetilde{C})= \mu_d(\alpha\,\widetilde{A}+\beta\,\widetilde{B})$,
the argument in the case of $A$ and $B$ yields the analogues of
(\ref{CAB+}) and  (\ref{CAB-}); namely,
\begin{eqnarray}
\label{tildeCAB+}
\widetilde{C}_+(t)&=&\alpha\,\widetilde{A}_+(t)+\beta\, \widetilde{B}_+(t),
\mbox{ \ and}\\
\nonumber
&&\mbox{$\widetilde{A}_+(t)$ and $\widetilde{B}_+(t)$ are $(d-1)$-balls with
$(d-1)$-volume $t$ for $t\in[t_+,1]$},\\
\label{tildeCAB-}
\widetilde{C}_-(t)&=&\alpha\,\widetilde{A}_-(t)+\beta\, \widetilde{B}_-(t),\mbox{ \ and}\\
\nonumber
&&\mbox{$\widetilde{A}_-(t)$ and $\widetilde{B}_-(t)$ are $(d-1)$-balls with
$(d-1)$-volume $t$ for $t\in[t_-,1]$}.
\end{eqnarray}

Let $\Pi$ be any two-dimensional linear subspace containing $u$. In particular,
$\Pi\cap\widetilde{A}=p_\Pi \widetilde{A}$,
$\Pi\cap\widetilde{B}=p_\Pi \widetilde{B}$ and $\Pi\cap\widetilde{C}=p_\Pi \widetilde{C}$,
and hence (\ref{projsum}) implies
\begin{equation}
\label{tildesection}
\Pi\cap\widetilde{C}=\alpha(\Pi\cap\widetilde{A})+\beta(\Pi\cap\widetilde{B}).
\end{equation}
We plan to apply Claim~\ref{planarsum} to $\Pi\cap\widetilde{A}$,
$\Pi\cap\widetilde{B}$ and $\Pi\cap\widetilde{C}$ with $l=\R u$.
Let $v\in S^{d-1}\cap u^\bot\cap \Pi$.
We observe that for $t_+<t\leq 1$, the radii of
$\widetilde{A}_+(t)$, $\widetilde{B}_+(t)$ and $\widetilde{C}_+(t)$
coincide by (\ref{tildeCAB+}), and
if $x\in \widetilde{A}_+(t)$, $y\in\widetilde{B}_+(t)$, $z\in \widetilde{C}_+(t)$
are relative boundary points with exterior normal $v$, then
there exists a common exterior unit normal vector to $\widetilde{A}$ at $x$ and
to $\widetilde{B}$ at $y$ by Claim~\ref{Minkowskisum}.
Combining this with the analogous properties of $\widetilde{A}_-(t)$, $\widetilde{B}_-(t)$
 and $\widetilde{C}_-(t)$ implies condition (i) of Claim~\ref{planarsum}.
In addition if $z_0\in \widetilde{C}_+(1)$ is
 relative boundary point with exterior normal $v$, then (\ref{AB1})
yields that $z_0+l$ intersects all of
$\widetilde{A}_+(1)$, $\widetilde{B}_+(1)$ and  $\widetilde{C}_+(1)$
in $\{z_0\}$. Therefore we may apply Claim~\ref{planarsum},
and deduce that $\Pi\cap\widetilde{A}=\Pi\cap\widetilde{B}$.
Therefore $\widetilde{A}=\widetilde{B}$, which in turn yields (\ref{abt}).

Next we claim that
\begin{equation}
\label{hAhB}
h_A(w)=h_ B(w) \mbox{ \ for $w\in S^{n-1}$}.
\end{equation}
We may assume that $w\neq \pm u$, and let $\Pi$ be
two-dimensional linear subspace spanned by $u$ and $w$.
Again let $v\in S^{d-1}\cap u^\bot\cap \Pi$.
We plan to apply Claim~\ref{planarsum} to $p_\Pi A$, $p_\Pi B$ and $p_\Pi C$ with $l=\R v$.
We deduce condition (i) by (\ref{CAB+}), (\ref{CAB-}) and (\ref{abt}),
and condition (ii) by (\ref{AB1}). Therefore
$p_\Pi A=p_\Pi B$, and hence $h_A(w)=h_{p_\Pi A}(w)=h_{p_\Pi B}(w)=h_ B(w)$.

Finally (\ref{hAhB}) yields that $A=B$.\\

\noindent{\bf Case 2 } Either $a_+(1)>a_-(1)$, or $b_+(1)>b_-(1)$.

We may assume that $a_+(1)-a_-(1)\geq b_+(1)-b_-(1)$, and hence
$a_+(1)>a_-(1)$. It follows that
$A_+(1)=A_-(1)+\lambda v$ for suitable $v\in S^{d-1}$
and $\lambda>0$. It follows by (\ref{Kstretch1}) that
\begin{equation}
\label{A1proj}
\pi_vA=\pi_vA_-(1).
\end{equation}

If $b_+(1)>b_-(1)$ then
$B_+(1)=B_-(1)+\tau w$ for $w\in S^{d-1}$
and $\tau>0$. If $b_+(1)=b_-(1)$, then we set $\tau=0$ and $w=v$,
and still have $B_+(1)=B_-(1)+\tau w$.
It follows from (\ref{CAB+}) and  (\ref{CAB-}) that
$C_+(1)=C_-(1)+\alpha\lambda v+\beta \tau w$.
We deduce by (\ref{Kstretch0}) and  (\ref{Cmax}) that
\begin{eqnarray*}
C_-(1)+[0,\alpha\lambda v]+[0,\beta \tau w]&=&
\alpha\left(A_-(1)+[0,\lambda v]\right)+
\beta\left(B_-(1)+[0, \tau w]\right)\\
&\subset&\{z\in C:\,[c_-(1)\leq\langle z,u\rangle\leq c_+(1)\}\\
&=&C_-(1)+[0,\alpha\lambda v+\beta \tau w].
\end{eqnarray*}
Since $\langle v,u\rangle>0$ and $\langle w,u\rangle>0$,
we conclude that $v=w$ also if $b_+(1)>b_-(1)$.

We deduce by (\ref{Kstretch0}) that
\begin{equation}
\label{ACRv}
A\subset A_-(1)+\R v\mbox{ and }C\subset C_-(1)+\R v,
\end{equation}
and claim that
\begin{equation}
\label{BRv}
B\subset B_-(1)+\R v.
\end{equation}
If $\tau>0$ then (\ref{BRv}) also follows from (\ref{Kstretch0}).
If $\tau =0$, then we should prove that $y_0+\R v$ is a supporting line to $B$
for any relative boundary point $y_0$ of $B_-(1)$.
Now $x_0=y_0$ is a relative boundary point of $A_-(1)=B_-(1)$ (see (\ref{AB1})),
hence $z_0=\alpha\,x_0+\beta\,y_0=y_0$
 is a relative boundary point of $C_-(1)=\alpha\,A_-(1)+\beta\,B_-(1)=B_-(1)$
(compare (\ref{CAB-})). Thus (\ref{ACRv}) yields that there exists
a supporting hyperplane $H$ containing $z_0+\R v$ at $z_0$ to $C$,
and in turn Claim~\ref{Minkowskisum}) (ii) implies that $H$ is
a supporting hyperplane at $y_0$ to $B$. We conclude (\ref{BRv}).

We define
\begin{eqnarray*}
A'&=&\bigcup_{t\in[0,1]}\left((A_+(t)-\lambda v)\cup A_-(t)\right)\\
B'&=&\bigcup_{t\in[0,1]}\left((B_+(t)-\tau v)\cup B_-(t)\right)\\
C'&=&\bigcup_{t\in[0,1]}\left((C_+(t)-\alpha\lambda v-\beta\tau v)\cup C_-(t)\right).
\end{eqnarray*}
We deduce by $C=A+B$,  (\ref{ACRv}) and (\ref{BRv}) that $C'=A'+B'$.
In addition
\begin{eqnarray*}
\mu_d(C')&=&\mu_d(C)-\mu_{d-1}(C_-(1))\cdot
\langle(\alpha\lambda +\beta\tau) v,u\rangle\\
&=&\alpha\mu_d(A)+\beta\mu_d(B)-
\alpha\mu_{d-1}(A_-(1))\cdot\langle\lambda v,u\rangle-
\beta\mu_{d-1}(A_-(1))\cdot\langle \tau v,u\rangle\\
&=&\alpha\mu_d(A')+\beta\mu_d(B').
\end{eqnarray*}
Since both $A'$ and $B'$ have a unique section parallel to $L$ of maximal $(d-1)$-dimensional
volume, we deduce by Case 1 that $A'=B'$.
We conclude  Theorem~\ref{Bonnesen-section-equa} by (\ref{A1proj}).
\proofbox

\section{Proof of Theorem~\ref{Bonnesen-proj-equa}}
\label{secproj}

In this section, we assume
\begin{equation}
\label{proj1}
M=\mu_{d-1}(\pi_u A)=\mu_{d-1}(\pi_u B)=N=1
\mbox{ \ and \ }\alpha+\beta=1.
\end{equation}

If the convex bodies $A'$ and $B'$ are homothetic, then (\ref{proj1}) yields that $A'$
and $B'$ are translates. If in addition $A$ and $B$ are obtained from $A'$
$B'$, respectively, by stretching along $u$, then readily
$$
\mu_d(\alpha A+\beta B)=\alpha\mu_d(A)+\beta\mu_d(B).
$$

For the reverse direction, first we explain how Theorem~\ref{Bonnesen-section} yields
 Theorem~\ref{Bonnesen-proj} {\it via} Steiner symmetrization.
Let $\widetilde{A}$, $\widetilde{B}$, and $\widetilde{C}$ be the Steiner symmetrials of
$A$, $B$ and $C=\alpha A+\beta B$. In particular
$\alpha \widetilde{A}+\beta \widetilde{B}\subset\widetilde{C}$ according to
Claim~\ref{steiner-symm}. We also observe that
$\pi_u A=u^\bot\cap \widetilde{A}$ and $\pi_u B=u^\bot\cap\widetilde{B}$
are sections of maximal $(d-1)$-measure of $\widetilde{A}$ and $\widetilde{B}$,
respectively, parallel to $L=u^\bot$. Therefore  Theorem~\ref{Bonnesen-section}
and the conditions (\ref{proj1}) yield
$$
\mu_d(\alpha A+\beta B)=\mu_d(\widetilde{C})\geq
\mu_d(\alpha \widetilde{A}+\beta \widetilde{B})\geq
\alpha\mu_d(\widetilde{A})+\beta\mu_d(\widetilde{B})=\alpha\mu_d(A)+\beta\mu_d(B).
$$

Next we assume that $\mu_d(\alpha A+\beta B)=\alpha\mu_d(A)+\beta\mu_d(B)$,
and hence
\begin{eqnarray}
\label{proj-equa-steiner}
\widetilde{C}&=&\alpha \widetilde{A}+\beta \widetilde{B}\\
\label{proj-steiner-equa}
\mu_d(\alpha \widetilde{A}+\beta \widetilde{B})&=&
\alpha\mu_d(\widetilde{A})+\beta\mu_d(\widetilde{B}).
\end{eqnarray}
Combining (\ref{proj-steiner-equa}) and  Theorem~\ref{Bonnesen-section-equa}
shows that there exist homothetic convex bodies $\widetilde{A}'$ and
$\widetilde{B}'$, and a $v\in S^{d-1}$
such that $\widetilde{A}$ and $\widetilde{B}$ are obtained from
$\widetilde{A}'$ and $\widetilde{B}'$, respectively, by stretching along $v$. Since
$\widetilde{A}$ and $\widetilde{B}$ are symmetric through $u^\bot$,
we deduce that $v=\pm u$. Therefore we may assume that
$\widetilde{A}'$ and $\widetilde{B}'$ are also symmetric through $u^\bot$.
We deduce by the conditions (\ref{proj1}) that actually $\widetilde{A}'$ and
$\widetilde{B}'$
are translates, therefore $\widetilde{A}'=\widetilde{B}'$ can be assumed.
Therefore there exists a non-negative convex function $\varphi$
on $\pi_uA=\pi_uB$, and $a,b\geq 0$, such that
\begin{eqnarray*}
\widetilde{A}'=\widetilde{B}'&=&
\{x+\lambda u:\,x\in\pi_uA\mbox{ and }-\varphi(x)\leq \lambda\leq \varphi(x)\}\\
\widetilde{A}&=&
\{x+\lambda u:\,x\in\pi_uA\mbox{ and }-\varphi(x)-a\leq \lambda\leq \varphi(x)+a\}\\
\widetilde{B}&=&
\{x+\lambda u:\,x\in\pi_uA\mbox{ and }-\varphi(x)-b\leq \lambda\leq \varphi(x)+b\}.
\end{eqnarray*}

We deduce by (\ref{steiner-func}) that there exist functions $\theta$
and $\psi$ on $\pi_uA$ such that
\begin{eqnarray*}
A&=&
\{x+\lambda u:\,x\in\pi_uA\mbox{ and }\theta(x)-\varphi(x)-a\leq \lambda\leq
\theta(x)+\varphi(x)+a\}\\
B&=&
\{x+\lambda u:\,x\in\pi_uA\mbox{ and }
\psi(x)-\varphi(x)-b\leq \lambda\leq \psi(x)+\varphi(x)+b\}.
\end{eqnarray*}
It follows that $\theta(x)+\varphi(x)+a$, $-(\theta(x)-\varphi(x)-a)$,
$\psi(x)+\varphi(x)+b$ and $\psi(x)-\varphi(x)-b$
are convex. Since convex functions on a compact set are Lipschitz,
both $\varphi$ and $\theta$  are almost everywhere differentiable on $\pi_uA$.
For each $x\in\pi_u{\rm int}\,A$,  there are parallel supporting hyperplanes
to $\widetilde{A}$ at $x+(\varphi(x)+a)u$, and
to $\widetilde{B}$ at $x+(\varphi(x)+b)u$, thus
 (\ref{proj-equa-steiner}) and Claim~\ref{steiner-symm} that
$$
(\theta(x)+\varphi(x)+a)'=(\psi(x)+\varphi(x)+b)'
\mbox{ \ for almost all $x\in\pi_uA$}.
$$
Therefore there exists some $\omega\in\R$ such that $\psi(x)=\theta(x)+\omega$
for $x\in\pi_uA$. By possibly interchanging the role of $A$ and $B$, we may assume that
$\omega\geq 0$.
In particular defining
$$
A'=B'=
\{x+\lambda u:\,x\in\pi_uA\mbox{ and }\theta(x)-\varphi(x)\leq \lambda\leq
\theta(x)+\varphi(x)\},
$$
both $A$ and $B$ are obtained from $A'=B'$ by stretching along $u$.
\proofbox

\end{document}